\documentclass[a4paper,10pt,reqno]{amsart}
\usepackage{verbatim}
\usepackage{amssymb}

\usepackage{enumerate}
\usepackage[active]{srcltx}
\numberwithin{equation}{section}
\makeatletter
\@namedef{subjclassname@2020}{%
  \textup{2020} Mathematics Subject Classification}
\makeatother

\usepackage{t1enc}
\usepackage[utf8x]{inputenc}

\newtheorem{theorem}{Theorem}[section]

\newtheorem{corollary}[theorem]{Corollary}

\theoremstyle{definition}

\theoremstyle{remark}

\newcommand{\mc}[1]{\mathcal{#1}}

\newcommand{\setm}{\setminus}

\newcommand{\subs}{\subset}

\def\<{\left\langle}
\def\>{\right\rangle}
\def\br#1;#2;{\bigl[ {#1} \bigr]^ {#2} }

\author[I. Juh\'asz]{Istv\'an Juh\'asz}
\address      { Alfr\'ed Rényi Institute of Mathematics}
\email{juhasz@renyi.hu}

\author[L. Soukup]{Lajos Soukup}
\thanks
  {
   }
\address
      { Alfr{\'e}d R{\'e}nyi Institute of Mathematics%
}
\email{soukup@renyi.hu}

\author[Z. Szentmikl\'ossy]{Zolt\'an Szentmikl\'ossy}
\address{E\"otv\"os University of Budapest}
\email{szentmiklossyz@gmail.com}


\title{On a problem of Angelo Bella}
\thanks{The research on and preparation of this paper was
supported by  NKFIH grant no.  K 129211.}
\date{\today}

\begin{document}

\begin{abstract}
The main result of this note is the following theorem.

{\em If $X$ is any Hausdorff space with $\kappa = \widehat{F}(X) \cdot \widehat{\mu}(X)$
then $L(X_{< \kappa}) \le \varrho(\kappa)$.}

Here $\widehat{F}(X)$ is the smallest cardinal $\varphi$ so that $|S| < \varphi$ for any set $S$
that is free in $X$ and $\widehat{\mu}(X)$ is the smallest cardinal $\mu$ so that, for every set 
$S$ that is free in $X$, any open cover of $\overline {S}$ has a subcover of size $< \mu$.
Moreover, $X_{< \kappa}$ is the $G_{< \kappa}$-modification of $X$ and 
$\varrho(\kappa) = \min \{\varrho : \varrho ^{< \kappa} = \varrho\}$.

As a corollary we obtain that if $X$ is a linearly Lindelöf regular space of countable tightness 
then $L(X_\delta) \le \mathfrak{c}$, provided that $ \mathfrak{c} = 2^{< \mathfrak{c}}$.
This yields a consistent affirmative answer to a question of Angelo Bella.
\end{abstract}

\maketitle

Some time ago, in an e-mail message, Angelo Bella asked us the following question:
Let $X$ be a linearly Lindelöf regular space of countable tightness. Is it true
then that in the $G_{\delta}$-modification $X_\delta$ of $X$ every set of regular
cardinality $> \mathfrak{c}$ has a complete accumulation point?

While, at the first sight, this question might seem rather ad hoc, it is not. Indeed,
we recall that $X$ being linearly Lindelöf is equivalent to every subset of $X$ of uncountable regular
cardinality having  a complete accumulation point. Thus it is natural to define the "linear Lindelöf
number" $\ell L(X)$ of a space $X$ as the smallest cardinal $\lambda$ such that every set in $X$ of regular
cardinality $ > \lambda$ has a complete accumulation point. Thus, using this notation, the question asks if
$\ell L(X) \cdot t(X) = \omega$ implies $\ell L(X_\delta) \le \mathfrak{c}$, or more generally, if
$$\ell L(X_\delta) \le 2^{\ell L(X) \cdot t(X)}$$
holds for any regular space $X$.

Now, this question is raised naturally in the line of research that asks about the determination of the values of various
cardinal functions on the $G_{\delta}$-modification $X_\delta$ in terms of ones defined on $X$.
A systematic study of this general question had been initiated in \cite{BS1}.

The aim of this note is to give a partial affirmative answer to Angelo Bella's question. For instance,
we shall actually prove that $\ell L(X) \cdot t(X) = \omega$ implies even $L(X_\delta) \le \mathfrak{c}$,
provided that $ \mathfrak{c} = 2^{< \mathfrak{c}}$.
This will be a consequence of a much more general result that we think is of independent interest.

In order to formulate this general result, we need to introduce some (partially new) notation and terminology
which are related to the free set number as it appeared in \cite{JSSz1}.
Otherwise, our notation and terminology follows \cite{Ju}.

As is well-known, a transfinite sequence $\<x_\alpha : \alpha < \eta\>$ of points of a
topological space $X$ is said to be a {\em free sequence} in $X$ if the closure of any
initial segment of it is disjoint from the closure of the corresponding final segment,
i.e. $$\overline{\{x_\alpha : \alpha < \beta\}} \cap \overline{\{x_\alpha : \alpha \ge \beta\}} = \emptyset$$
for each $\beta < \eta$.
Now, we call a {\em subset} $S$ of a space $X$ {\em free} in $X$ if it admits a well-ordering,
or equivalently an indexing by ordinals,
that turns it into a free sequence in $X$.
In other words, free sets in $X$ are just the ranges of free sequences in $X$. We shall use $\mc{F}(X)$ to denote the collection of
all free subsets in $X$.

While in \cite{JSSz1} the main object of study was the free set number $$F(X) = \sup \{|S| : S \in \mc{F}(X)\},$$
here we shall use its "hat version" $\widehat{F}(X)$ that is defined as the smallest cardinal $\kappa$
such that $|S| < \kappa$ for all $S \in \mc{F}(X)$.

Similarly, instead of the auxiliary cardinal function $${\mu}(X)=\sup\{L(\overline {S}):S\in \mc F(X)\}$$
from \cite{JSSz1} we shall use its "hat version" $\widehat{\mu}(X)$ that is defined as the smallest cardinal $\kappa$
such that for every $S \in \mc{F}(X)$ and every open cover $\mathcal{U}$ of $\overline {S}$ there is $\mathcal{V} \subs \mathcal{U}$
with $|\mathcal{V}| < \kappa$ that also covers $\overline {S}$.

Finally, for any given uncountable cardinal $\kappa$, instead of the $G_{\delta}$-modification $X_\delta$ of a space $X$
we are going to consider the $(<\kappa)$-modification $X_{< \kappa}$,
which is the topology on $X$ generated by all $G_{< \kappa}$-sets in $X$, i.e. sets obtainable as intersections of fewer than $\kappa$
open subsets of $X$.

We are now ready to formulate our main result.

\begin{theorem}\label{tm:main}
Let $X$ be any Hausdorff space and set $\kappa = \widehat{F}(X) \cdot \widehat{\mu}(X)$. Then for every cardinal
$\varrho$ satisfying $\varrho ^{< \kappa} = \varrho$ we have $L(X_{< \kappa}) \le \varrho$.
\end{theorem}

\begin{proof}
Our proof is indirect, so we assume that $L(X_{< \kappa}) > \varrho$. This means that to every point $x \in X$
we may associate a family $\mathcal{U}_x$ of open sets in $X$ with $ |\mathcal{U}_x| < \kappa$ so that $x \in H_x = \bigcap \mathcal{U}_x$,
moreover no subfamily of $\mathcal{H} = \{H_x : x \in X\}$ of cardinality $\le \varrho$ covers $X$.

Next we consider an elementary submodel
$M$ of $H(\vartheta)$ for a large enough regular cardinal $\vartheta$ such that $|M| = \varrho$, $\,M$ is
$(<\kappa)$-closed and $\varrho + 1 \subs M$, moreover both $X$ and the map $x \mapsto \mathcal{U}_x$ belong to $M$.
There is such an $M$ because $\varrho ^{< \kappa} = \varrho$.

Then for $\mathcal{K} = \mathcal{H} \cap M$ we have $|\mathcal{K}| \le \varrho$, hence we may pick and fix a point $p \in X \setm \bigcup \mathcal{K}$.
We then define $\mathcal{W} = \{U \in M \cap \tau(X) : p \notin U\}$. Clerly, $\mathcal{W}$ is closed under unions
of size $< \kappa$ because $\,M$ is $(<\kappa)$-closed.

Next we want to show that for every $S \in \mc{F}(X)$ with $S \subs M$ there is $W \in \mathcal{W}$ so that $\overline {S} \subs W$;
this is really the crux of our proof.
To see this, we first note that $S \in M$ because $\,M$ is $(<\kappa)$-closed, consequently we have $\overline{S} \in M$ as well.
We then observe that, as $X$ is Hausdorff, $\overline{S}$ has a pseudo-base $\mathcal{B} \subs \tau(\overline{S})$ of size $\le 2^{|S|}$,
moreover $2^{|S|} \le \varrho$ because $|S| < \widehat{F}(X) \le \kappa$ and $\varrho ^{< \kappa} = \varrho$.

We next show that $\overline{S}$ admits a network $\mathcal{N}$ of size $\le \varrho$. Indeed, fix $x \in U \in \tau(\overline{S})$.
Then for every point $y \in \overline{S} \setm U$ we may pick $B_y \in \mathcal{B}$ such that $y \in B_y$ and $x \notin B_y$.
Now, it follows from the definition of  $\widehat{\mu}(X)$ that there is a subset $Y \subs \overline{S} \setm U$ with
$|Y| < \widehat{\mu}(X) \le \kappa$ such that $\overline{S} \setm U \subs \bigcup \{B_y : y \in Y\}$. But this clearly
implies that if we let $\mathcal{C} = \{\overline{S} \setm B : B \in \mathcal{B}\}$ then
$$\mathcal{N} = \{\cap \mathcal{D} : \mathcal{D} \in \big[\mathcal{C}\big]^{< \kappa}\}$$
is indeed a network for $\overline{S}$ with  $|\mathcal{N}| \le \varrho ^{< \kappa} = \varrho$.

But it is clear that then all intersections of fewer than $\kappa$ many members of $\mathcal{N}$ form
a network for the $G_{< \kappa}$-modification $(\overline{S})_{< \kappa}$ of $\overline{S}$,
consequently we have $L\big((\overline{S})_{< \kappa}\big) \le \varrho$.

It follows from this that the cover $\{H_x : x \in \overline{S}\} \subs \tau(X_{< \kappa})$ of $\overline{S}$
has a subcover of cardinality $\le \varrho$, i.e. there is some $T \subs \overline{S}$ with $|T| \le \varrho$
such that $\overline{S} \subs \bigcup \{H_x : x \in T\}$. But because $\overline{S} \in M$ and
the map $x \mapsto H_x$ belongs to $M$, by elementarity we may assume that $T \in M$. Note that then,
as $|T| \le \varrho$, we have $T \subs M$ as well.

Now, this implies that $H_x \in \mathcal{K}$ for every $x \in T$, hence we have $p \notin H_x$.
This, in turn, implies that for every $x \in T$ there is $U_x \in \mathcal{U}_x$ such that $p \notin U_x$.
But for any $x \in T$ we have $x \in M$, hence $\mathcal{U}_x \in M$ and so $U_x \in M$,
and this clearly implies $U_x \in \mathcal{W}$ as well.

Thus we may apply again our definition of  $\widehat{\mu}(X)$ to the open cover $\{U_x : x \in T\} \subs \tau(X)$ of
$\overline{S}$ to find $Y \subs T$ with $|Y| < \widehat{\mu}(X) \le \kappa$ so that  $\overline{S} \subs W = \bigcup \{U_x : x \in Y\}$.
But we have already seen that $\mathcal{W}$ is closed under unions
of size $< \kappa$, hence $W \in \mathcal{W}$, completing the proof of our claim.

Finally, it is obvious from elementarity that for every member $W \in \mathcal{W}$ we have $X \cap M \setm W \ne \emptyset$.
Consequently, the triple $\<X \cap M,\,\widehat{F}(X),\,\mathcal{W}\>$ satisfies all three conditions of Lemma 2.1
of \cite{JSSz1}, hence there is a subset $S \subs X \cap M$ of cardinality $\widehat{F}(X)$ that is free in $X$.
This, however, is a blatant contradiction that completes the proof of our theorem.
\end{proof}

For every uncountable cardinal $\kappa$ we shall denote by $\varrho(\kappa)$ the smallest cardinal $\varrho$
such that $\varrho ^{< \kappa} = \varrho$. Clearly, then our above theorem could also be phrased as follows:

{\em For any Hausdorff space $X$ with $\kappa = \widehat{F}(X) \cdot \widehat{\mu}(X)$
we have $L(X_{< \kappa}) \le \varrho(\kappa)$.}

\smallskip

So, let us now examine how $\varrho(\kappa)$ can be calculated. First, it is obvious
that $2^{< \kappa} \le \varrho(\kappa) \le 2^\kappa$. It is also easy to see that if
$\kappa$ is regular then $\varrho(\kappa) = 2^{< \kappa}$. Indeed, this immediately follows
from $cf(2^{< \kappa}) \ge \kappa$.

If $\kappa$ is singular and there is a $\mu < \kappa$ with $2^\mu  = 2^{< \kappa}$ then again
we obviously have $\varrho(\kappa) = 2^{< \kappa}$.

So, we are left with the case in which $\kappa$ is singular and for all $\mu < \kappa$ we have $2^\mu  < 2^{< \kappa}$.
We claim that in this case $\varrho(\kappa) = 2^\kappa$.

Indeed, in this case there is a strictly increasing sequence of cardinals $\<\kappa_i : i < cf(\kappa)\>$ converging to $\kappa$ 
such that the (increasing) sequence $\<\lambda_i = 2^{\kappa_i} : i < cf(\kappa)\>$ converges to $\lambda = 2^{< \kappa}$.
In particular, we then have $cf(\lambda) = cf(\kappa)$.
It is well known, however, that then we have $\lambda^{cf(\kappa)} = \prod \{\lambda_i : i < cf(\kappa)\}$.
But, on the other hand, we have 
$$\prod \{\lambda_i : i < cf(\kappa)\} = \prod \{2^{\kappa_i} : i < cf(\kappa)\} = 2^{\sum \{\kappa_i\, :\, i < cf(\kappa)\}} = 2^\kappa.$$
Putting all these things together we get $$2^\kappa = \lambda^{cf(\kappa)} \le \varrho(\kappa)^{cf(\kappa)} = \varrho(\kappa) \le 2^\kappa,$$
hence, indeed,  $\varrho(\kappa) = 2^\kappa$.

\medskip

Now we turn to discussing how our Theorem \ref{tm:main} yields a partial affirmative answer to Bella's above question.
First of all, we note that for any linearly Lindelöf space $X$ of countable tightness we have $F(X) = \omega$ because then
a free sequence of length $\omega_1$ in $X$ cannot have a complete accumulation point.

Secondly, if $X$ is also regular then we claim that we also have $\widehat{\mu}(X) \le \mathfrak{c}$.
Indeed, since every $S \in \mathcal{F}(X)$ is countable, we have $L(\overline{S}) \le w(\overline{S}) \le \mathfrak{c}$.
But because $X$ and hence $\overline{S}$ is linearly Lindelöf, every open cover of $\overline{S}$ has a subcover
whose cardinality has countable cofinality. Consequently, since $cf(\mathfrak{c}) > \omega$, we actually have
$$\widehat{\mu}(X) \le \sup \{\lambda^+ : \lambda < \mathfrak{c} \text{ and } cf(\lambda) = \omega\} \le \mathfrak{c}.$$
This, of course, immediately implies $\widehat{F}(X) \cdot \widehat{\mu}(X) \le \mathfrak{c}$ and thus the promised 
consistent affirmative answer to Angelo Bella's question.

\begin{corollary}\label{co:w}
If $X$ is any linearly Lindelöf regular space of countable tightness then $\ell L(X_\delta) \le L(X_\delta) \le \varrho(\mathfrak{c})$.
In particular, if $\mathfrak{c} = 2^{< \mathfrak{c}}$ then $\ell L(X_\delta) \le L(X_\delta) \le \mathfrak{c}$.
\end{corollary}

It remains a wide open question whether the affirmative answer to Angelo Bella's question is provable in ZFC or not.
In fact, we don't even know if the first sentence of Corollary \ref{co:w} remains valid in ZFC when $\varrho(\mathfrak{c})$ is replaced by $2^{< \mathfrak{c}}$.
However, it is clear that the argument that lead us to Corollary \ref{co:w} can be generalized to obtain
the following result. 

\begin{corollary}\label{co:gen}
For any regular space $X$ we have $$\ell L(X_\delta) \le L(X_\delta) \le \varrho(\ell L(X) \cdot t(X)).$$
\end{corollary}

\end{document}